%
%
\input amstex.tex
\input amsppt.sty
\input xy
\xyoption{all}
\documentstyle{amsppt}
\NoBlackBoxes

\def\LL{\Bbb L}
\def\SS{\Cal S}


\topmatter
\title 
Controlled homotopy equivalences and structure sets of manifolds
\endtitle

\author 
Friedrich Hegenbarth and Du\v san Repov\v s
\endauthor

\leftheadtext{F. Hegenbarth and D. Repov\v s}
\rightheadtext{Controlled homotopy equivalences and structure sets}

\address
Department of Mathematics,
University of Milano, 
Via C. Saldini 50, 
Milano, Italy 02130.
\endaddress
\email
friedrich.hegenbarth\@mat.unimi.it 
\endemail

\address
Faculty of Education, 
and Faculty of Mathematics and Physics,
University of Ljubljana, Kardeljeva pl. 16, 
Ljubljana, Slovenia 1000.
\endaddress
\email
dusan.repovs\@guest.arnes.si 
\endemail




\abstract
For a closed topological $n$--manifold $K$ and a
map $p:K\to B$ inducing an isomorphism $\pi_1(K)\to\pi_1(B)$, there is a canonicaly defined morphism $b:H_{n+1}(B,K,\LL)\to \SS (K)$, where $\LL$ is the periodic simply-connected surgery spectrum and $\SS (K)$ is the topological structure set. We construct a refinement $a:H_{n+1}^{+}(B,K,\LL )\to \SS_{\varepsilon ,\delta }(K)$ in the case when $p$ is $UV^1$, and we show that $a$ is bijective  if $B$ is a finite-dimensional compact metric ANR. 
Here, $H_{n+1}^{+}(B,K,\LL )\subset H_{n+1}(B,K,\LL )$, and $\SS_{\varepsilon ,\delta }(K)$ is the controlled structure set. We show that the Pedersen-Quinn-Ranicki
controlled surgery sequence  is equivalent to the exact $\LL$-homology sequence of the map $p:K \to B$, i.e. that
$$H_{n+1}(B,\LL)\to H_{n+1}^{+}(B,K,\LL )\to H_n(K,\LL^{+})\to H_n(B,\LL ), \ 
\LL^{+}\to \LL,$$ is the connected covering spectrum of $\LL$. By taking for $B$ various stages of the Postnikov tower of $K$, one obtains an interesting filtration of the controlled structure set.

\endabstract

\endtopmatter

\document

\head {$\S$ 1. Introduction}\endhead
Let $K$ denote a closed topological manifold of dimension $n$. As usual, $\SS(K)$ denotes the topological structure set of $K$.  Elements of
$\SS(K)$ are equivalence classes of pairs $(M,h)$, where $M$ is a closed topological $n$-manifold and $h:M\to K$ is 
a simple homotopy  equivalence. The pairs $(M,h),(M_1,h_1)$ are equivalent if there is a homeomorphism  $\varphi :M_1\to M$ such that $h\circ \varphi$ is homotopic to $h_1$. $\SS(K)$ has a group structure. If $B$ is a finite-dimensional compact metric ANR and $p:K\to B$ a continuous $UV^1$ map, the controlled structure set $\SS_{\varepsilon ,\delta }(K\to B)$  (cf. \cite{PQR}) of $\delta$-homotopy equivalences can be defined (see
precise definitions below). There is an obvious forgetful map $\SS_{\varepsilon ,\delta }(K\to B)\to \SS(K)$.
An element of
the image will be  called a controlled homotopy equivalence
$h:M\to K$ with respect to $p$. Varying $p:K\to B$, we get the subset of controlled homotopy equivalences
in $\SS (K)$.

Dranishnikov and Ferry studied the
subset (indeed a subgroup) $\SS^{CE}(K)\subset \SS(K)$ of elements realized by cell-like maps (cf. \cite{DrFe}). They observed that these elements belong to the image of a natural map
$$
b:H_{n+1}(E_2(K),K;\LL )\to \SS(K)$$
localized away from $2$.

Here $\LL$ is the simply-connected $4$-periodic surgery spectrum, 
and $E_2(K)$ is the second stage of
the Postnikov tower of $K$.
The map $b$ can be defined using the diagram on
p. 207 of \cite{Ran}. 
The main ingredient is the $\pi -\pi$--theorem. 
The morphism $b$ can be defined for any map $p:K\to B$
satisfying $p_*:\pi_1(K)\mathop\to\limits_{\cong} \pi_1(B)$, 
we shall generically call it  $b$.

One may ask of what type are homotopy equivalences of images of other Postnikov stages.
In this paper we shall give an answer in the
case of the first stage, i.e. $B\pi$, where $\pi=\pi_1(K)$. 
As a consequence of the main theorem we shall get the following result:

\proclaim{Theorem 1}
Let $n\geq 5$. 
Then the image of the map
$$
b:H_{n+1}(B\pi , K, \LL )\to \SS(K)
$$
contains the controlled homotopy equivalences.
\endproclaim

In order to be able to state the main result
we shall need some more notations. 
Recall that $\LL_0=\Bbb Z \times G/TOP$. Let $\LL^+\to \LL$ be the simply-connected spectrum  covering $\LL$. 
In particular, $\LL_0^+=G/TOP$. For any pair $(Z,Y)$ there is the following 
exact sequence
(cf. \cite{Ran, p. 153})
$$(*)\dots \to H_{m+2}(Z,Y,\Bbb Z)\to H_{m+1}(Z,Y,\LL^+)\to H_{m+1}(Z,Y,\LL )\to H_{m+1}(Z,Y,\Bbb Z)\to \dots $$

In particular,
for $(K,\emptyset)$, where $K$ is an oriented closed topological $n$-manifold, we have an injection
$$
H_n(K,\LL^+)\to H_n(K,\LL).
$$

Let us denote the inverse image of $H_n(K,\LL^+)$ under $\partial_*:H_{n+1}(B,K,\LL )\to 
H_n(K,\LL )$
by
$H_{n+1}^+(B,K,\LL )\subset H_{n+1}(B,K,\LL )$.
Here $B$ is a space containing $K$, or more generally a map $K\to B$.

\proclaim{Theorem 2}
Suppose that
$p:K\to B$ is a $UV^1$-map into a compact metric ANR-space, with $K$ as above and
$n\geq 5$. 
Then there is a map
$$
a:H_{n+1}^+(B,K,\LL )\to \SS_{\varepsilon,\delta}(K\to B),
$$
where $\varepsilon,\delta$ are appropriately choosen (see below). 
The map $a$ fits into the following
diagram
$$
\xymatrix @R=.1cm{
& H_{n+1}^+(B,K,\LL ) \ar[r]^{a} \ar[dd] & \SS_{\varepsilon,\delta}(K\to B) \ar[dd] \\
(**) & & \\
& H_{n+1}(B,K,\LL ) \ar[r]^{b} & \SS(K).
}
$$
\endproclaim

The natural map $H_{n+1}(B,\LL )\to H_{n+1}(B,K,\LL )$ factors over $H_{n+1}^+(B,K,\LL )$, and one can easily see that the sequence
$$
...\to H_{n+1}(B,\LL )\to H_{n+1}^+(B,K,\LL )\to H_{n}(K,\LL^+)\to H_{n}(B,\LL )\to ...
$$
is exact.

Recall that by the Poincar\'{e} duality $H_{n}(K,\LL^+)\cong [K,G/TOP]$.

\proclaim{Theorem 3}
Let $K$ be a closed oriented topological $n$-manifold, 
$n\geq 5$, and $p:K\to B$ a $UV^1$-map into a finite-dimensional 
compact metric ANR. Then there exists a commutative diagram of exact sequences 
$$
\xymatrix @C=.1cm @R=.8cm{
H_{n+1}(K,\LL^+) \ar[d]^{\cong} \to & H_{n+1}(B,\LL) \ar@{=}[d] & \to & H_{n+1}^+(B,K,\LL) \ar[d]^-{a} & \to & H_n(K,\LL^+) \ar[d]^{\cong} & \to & H_n(B,\LL) \ar@{=}[d] \\
[SX,G/TOP] \to  & H_n(B,\LL) & \to & \SS_{\varepsilon,\delta}(K\to B) & \to & [K,G/TOP] & \to & H_n(B,\LL).}
$$
In particular, $a$ is bijective.
\endproclaim

Here the lower exact sequence is the controlled surgery sequence from \cite{PQR} and
$\varepsilon$ and $\delta$ are to be appropriately choosen: 
There is an $\varepsilon_0 >0$, depending on $B$ and $n$, 
such that for any $0<\varepsilon < \varepsilon_0$ 
there is $\delta > 0$ such that the lower sequence is exact.

Theorem 2 and Theorem 3 will be proved in $\S 5$ and $\S 6$, respectively.
Theorem 1 follows from Theorem 2, since any $UV^1$-map $K\to B$
composed with the canonical map $B\to B\pi$ determines a 
map $(B,K)\to (B\pi,K)$. Since $K$ is a manifold of dimension
$\geq 5$ we can assume that $K\to B\pi$ is $UV^1$
(cf. \cite{BFMW}). Hence we get an induced morphism
$$
H_{n+1}(B,K,\LL)\to H_{n+1}(B\pi,K,\LL).
$$

On the other hand, any $[y]\in H_{n+1}(B\pi,K,\LL)$ has a 
compact support $(B,K)\subset (B\pi,K)$ (cf. \cite{Ran, 
$\S 12$}), and we may assume $K\to B$ is $UV^1$.

To define the map $a$ we recall the geometric $\LL$-spectrum, 
and represent elements $[x]\in H_{n}(K,\LL)$ and $[y]\in H_{n+1}(B,K,\LL)$ 
in terms of normal degree $1$ maps. We shall
follow \cite{Qu1}, \cite{Nic} and \cite{Ran}. 
One reason to recall these is to stabilize the
notations used  in the proofs. 
Another reason is that \cite{Ran} used the algebraic $\LL$-spectrum 
defined by quadratic Poincar\'{e} chain complexes which led to representing
$[x]$ by a normal degree $1$ map between closed manifolds (cf.
Corollary 18.6(I)
in \cite{Ran}). The gluing construction using the geometric $\LL$-spectrum
gives surgery problems with homotopy equivalences on boundaries. Explicit
examples have recently been given (cf. \cite{Ham, Example 5.4}).

\head {$\S$ 2. The periodic simply connected $\LL$-spectrum}\endhead

We briefly recall the semi-simplicial surgery spectrum introduced in \cite{Qu1}
(cf. also \cite{Nic}), denoted by $\LL=\{\LL_g|g\in \Bbb Z\}$.  It is an $\Omega$-spectrum,
and each $\LL_g$ satisfies the Kan condition. Its homotopy groups are the Wall groups
of the trivial fundamental group, more precisely, $\pi_j(\LL_g)=L_{j+g}(\{1\})$, and 
$\LL_0=\Bbb Z\times G/TOP$.

Let $\LL_g(j)$ denote the $j$-simplices of $\LL_g$. 
An element $x$ in $\LL_g(j)$ is 
a degree $1$ normal map of $(j+3)$-ads, which we
shall
shortly  denote by
$$
x=\{(f,b):M\to X\},
$$
where $(M, \partial _0 M,...,\partial _j M, \partial _{j+1}M)$ 
(resp. $(X,\partial _0 X,...,\partial _j X,\partial _{j+1} X)$)
are the $(j+3)$-ads $M$
(resp. $X$).
It is required that $f|_{\partial_{j+1}M}$ is a (simple) homotopy equivalence. $M$ and $X$
are topological manifolds of dimension $j+g$ with  boundaries 
$\partial X=\mathop\cup\limits_{i=0}^{j+1}\partial_iX$ and
$\partial M=\mathop\cup\limits_{i=0}^{j+1}\partial_iM$. 
To $x$ belongs a reference map from $X$ to the $(j+3)$-ad $\Delta^j$,
i.e. a map $(X,\partial _0 X,...,\partial _j X,\partial _{j+1} X)\to  (\Delta^j,\partial_0\Delta^j,...\partial_j\Delta^j,\Delta^j)$. 
We emphasize that the condition that
$f|_{\partial_{j+1}M}:\partial_{j+1}M\to\partial_{j+1}X$ 
is a homotopy equivalence is necessary to indentify the homotopy groups of $\LL$ 
with the Wall groups.

Cartesian product with $\Bbb C P^2$ defines a semi-simplicial map $\LL_g\to\LL_{g+4}$ 
which is a homotopy equivalence, called the periodicity map.


\head {$\S$ 3. $\LL$-homology classes}\endhead
We shall
follow \cite{Ran, $\S 12$}, to describe elements $[x]\in H_n(K, \LL)$, 
where $K$ is a finite ordered simplicial complex. The relative 
case will be given in $\S 4$.

The class $[x]$ is represented by a "cycle" $x$, which roughly is a 
semi-simplicial map from $K$ to $\LL_g$ for some $g$.  More precisely,
one embedds $K$ into $\partial\Delta^{m+1}$, $\Delta^{m+1}$ the standard
$(m+1)$-simplex 
with vertices denoted $\{0,1,...m+1\}$. The embedding is order-preserving and 
one may assume $m+1=$ number of vertices of $K$. For any $\sigma\in\partial\Delta^{m+1}$, 
let $D(\sigma,\partial\Delta^{m+1})$ be the dual cell, and $\Sigma^m$ 
the dual cell-complex decomposition of $\partial\Delta^{m+1}$. 
Note that $D(\sigma,\partial\Delta^{m+1})$ is an $(m-|\sigma|)$-ball with boundaries $\partial_iD(\sigma,\partial\Delta^{m+1})=D(\delta_i\sigma,\partial\Delta^{m+1})$, 
where $\sigma=\{j_0,...,j_k\}\subset\{0,...,m+1\}$, and $\delta_i\sigma=\sigma\cup\{\lambda_i\}$ 
with $\{0,...,m+1\}\setminus\{j_0,...,j_h\}= \{\lambda_0,...,\lambda _{m-k}\}$ (written in that order).
Hence $\{D(\sigma,\partial\Delta^{m+1}),\partial_0D(\sigma,\partial\Delta^{m+1}),..., \partial_{m-|\sigma|}D(\sigma,\partial\Delta^{m+1})\}$ is an $(m-|\sigma|+2)$-ball-ad.

Define $\overline{K}=\{D(\sigma,\partial\Delta^{m+1})|\sigma\in\partial\Delta^{m+1}\setminus K\}$
and call it the supplement of $K$.
The cycle $x$ can be considered as a semi-simplicial map
$$
x:(\Sigma^m,\overline{K})\to (\LL_{n-m},\emptyset) 
$$
i.e.
$x(D(\sigma,\partial\Delta^{m+1})\in \LL_{n-m}(m-|\sigma|)$ is a surgery problem of
$(m-|\sigma|+3)$-ads denoted simply
$$
x=\{x_{\sigma}=(f_{\sigma},b_{\sigma}):M_{\sigma}\to X_{\sigma}|\sigma\in\partial\Delta^{m+1}\}$$ 
between manifold $(m-|\sigma|+3)$-ads of dimension $n-|\sigma|$ together with a reference
map of $(n-|\sigma|+3)$-ads
$$
(X_\sigma,\partial_0X_\sigma,...,\partial_{m-|\sigma|}X_\sigma,\partial_{m-|\sigma|+1}X_\sigma) $$
$$\to (D(\sigma,\partial\Delta^{m+1}),\partial_0D(\sigma,\partial\Delta^{m+1}),...,\partial_{m-|\sigma|}D(\sigma,\partial\Delta^{m+1}),D(\sigma,\partial\Delta^{m+1})).
$$
Moreover , $f_{\sigma}$ restricted to $\partial_{m-|\sigma|+1}M_\sigma\to \partial_{m-|\sigma|+1}X_\sigma$ 
is a homotopy equivalence. Note that $X_\sigma\not=\emptyset\Rightarrow\sigma\in K$.

If $x'=\{x'_\sigma|\sigma\in\partial\Delta^{m+1}\}$ is a cycle defining the same class as $x$,
there is a
semi-simplicial map
$$
y=\Delta^1\times(\Sigma^m,\overline{K})\to (\LL_{n-m},\emptyset)
$$
such that $y(\partial_0\Delta^1\times\sigma)=x(\sigma)$, $y(\partial_1\Delta^1\times\sigma)=x'(\sigma)$.

Next we assemble the pieces $\{X_\sigma\}$ resp. $\{M_\sigma\}$, i.e. we build colimits. 
For this we use results from \cite{LaMc}  which are elegant reformulation of \cite{BRS}.

Let $A_{STOP}$ be the graded category with objects $n$-dimensional oriented topological 
manifolds with boundaries. The morphisms are orientation preserving inclusions 
(if dimension preserving) or inclusions with the
image in boundaries (if  dimension increasing) --
we refer to Example 3.5 from \cite{LaMc}.

Given the cycle $x$ as above we obtain two functors $(\Sigma^m,\overline{K})\to (A_{STOP},\emptyset)$ 
of degree $m-n$, namely we associate to $x(\sigma)$ either the target or the domain: i.e.
$$ D(\sigma,\partial\Delta^{m+1})\to X_\sigma \text{ and } D(\sigma,\partial\Delta^{m+1})\to M_\sigma .
$$
Let $X$ resp. $M$ denote its colimits, and denote by $(f,b):M\to X$ the resulting map. 
By Proposition 6.6 of \cite{LaMc},  $M$
and
$X$ are $n$-manifolds with boundaries $\partial M$, $\partial X$,
respectively. In fact,  $\partial M$ and $\partial X$ are the colimits of  the restricted 
functors to the boundary components $\partial_{m-|\sigma|+1}M_\sigma$
resp. $\partial_{m-|\sigma|+1}X_\sigma$.

There are no more boundary components. To see this we may consider a (simplicial) 
neighborhood $N(K)$ of $K\subset\partial\Delta^{m+1}$. 
Then
$X_\sigma=\emptyset$
for all $\sigma\not\in K$, 
in
particular $X_\sigma=\emptyset$ if $\sigma\in\partial N(K)$.

Together with the colimit of the reference maps one has obtained a normal degree $1$ 
maps $(f,b):M\to X$ with a reference map  $\varphi:X\to N(K)\to K$, where $N(K)\to K$ is a retraction.

\remark{Remarks}
\roster
\item "(I)"
The restriction of $f$ to $\partial M$ is a homotopy equivalence. Indeed,
it is a $\delta$-homotopy equivalence. This can be seen applying the gluing
construction
(cf. e.g. \cite{Hat, Ch. 4.G})
inductively to the
family $\{f_\sigma|...:\partial_{m-|\sigma|+1}M_\sigma\to \partial_{m-|\sigma|+1}X_\sigma|\sigma\in K\}$.
Here $\delta$ can be choosen as small as necessary by taking
a sufficiently fine subdivision.
\item "(II)"
A careful construction shows that $\varphi$ is transverse to all $D(\sigma,K)$, the 
duals of $\sigma$ in $K$, and that
$$
\varphi^{-1}(D(\sigma,K))=X_\sigma.$$
Conversely, any simplicial map from a manifold $\varphi:X\to K$ is transversal to 
dual cells (\cite{Coh}). Hence any normal degree $1$ map $M\to X$ defines a class $[x]\in H_n(K,\LL)$.
\item "(III)"
If $K$ is not a finite complex
then
one can choose for any $[x]\in H_n(K,\LL)$ a finite
complex $J\subset K$, so that $[x]$ belongs to the image
of $H_n(J,\LL)\to H_n(K,\LL)$.
\item "(IV)"
Instead of Prop. 6.6 from
\cite{LaMc} one can invoke the gluing procedure applied 
in \cite{Nic, Ch. 3}. However, the notion of degree decreasing functors from \cite{LaMc}
is more appropriate here.
\endroster
\endremark

\head {$\S$ 4. Relative $\LL$-homology classes}\endhead
Here we consider ordered simplicial complexes $B$, $K$ and a simplicial map $p:K\to B$. 
We assume $B$, $K$ to be finite complexes, and as in  $\S 3$ let $B\subset\partial\Delta^{m+1}$ 
be an order preserving embedding. By the
simplicial mapping cylinder construction we 
substitute $p$ by an inclusion $K\subset B$. An element $[y]\in H_{n+1}(B,K,\LL)$ 
is given by a semi-simplicial map
$$
y:(\overline{K},\overline{B})\to (\LL_{n+1-m},\emptyset)\text{, i.e.}
$$
by a family $\{y(\tau)\in\LL_{n+1-m}(m-|\tau|)\}$.

Each $y(\tau)$ is a surgery problem $(f_\tau, b_\tau):W_\tau\to V_\tau$ of manifold $(m-|\tau|+3)$-ads with $dimV_\tau=dimW_\tau=n+1-|\tau|$, such that $f_\tau$ restricts to a homotopy equivalence
$$
\partial_{m-|\tau|+1}W_\tau\to \partial_{m-|\tau|+1}V_\tau.
$$

If is convenient to give the family $\{y(\tau)\}$ a "cobordism" interpretation: 
$y$ can be written as semi-simplicial map
$$
y:(\Delta^1\times\Sigma^m,\Delta^1\times\overline{B}\cup\partial_0\Delta^1\times\overline{K}\cup\partial_1\Delta^1\times\Sigma^m)\to(\LL_{n-m},\emptyset)
$$
(cf. \cite{Ran, p. 128}).

The connecting homomorphism $\partial_*:H_{n+1}(B,K,\LL)\to H_n(K, \LL)$ 
sends $[y]$ to $\partial_*[y]=[x]$, where 
$$
x:(\partial_0\Delta^1\times\Sigma^m,\partial_0\Delta^1\times\overline{K})\to(\LL_{n-m},\emptyset)
$$
is the restriction of $y$.

If we write $\Delta^1\times\sigma$ for $\Delta^1\times D(\sigma,\partial\Delta^{m+1})\subset\Delta^1\times\Sigma^m$, considered as $(m-|\sigma|+1)-ball$, we have $y(\Delta^1\times\sigma)\in\LL_{n-m}(m-|\sigma|+1)$, i.e. it is a surgery problem
$$
(f_{\Delta^1\times\sigma},b_{\Delta^1\times\sigma}):W^{n+1-|\sigma|}_{\Delta^1\times\sigma}\to V^{n+1-|\sigma|}_{\Delta^1\times\sigma} 
$$
between $(m-|\sigma|+5)$-ads, such that $W_{\Delta^1\times\sigma}$, $V_{\Delta^1\times\sigma}$ are $(n+1-|\sigma|)$-dimensional manifolds and $f_{\Delta^1\times\sigma}$ restricts to a homotopy equivalence $$
\partial_{m-|\sigma|+3}W_{\Delta^1\times\sigma}\to \partial_{m-|\sigma|+3}V_{\Delta^1\times\sigma}.
$$

Note that $\partial_0V_{\Delta^1\times\sigma}=V_{\partial_0\Delta^1\times\sigma}$, $\partial_0W_{\Delta^1\times\sigma}=W_{\partial_0\Delta^1\times\sigma}$.
Furthermore, we have
\roster
\item "(I)"
$V_{\Delta^1\times\sigma}\not=\emptyset\Rightarrow \sigma\in B$;
\item "(II)"
$\partial_0V_{\Delta^1\times\sigma}\not=\emptyset\Rightarrow \sigma\in K$; and
\item "(III)"
$\partial_1V_{\Delta^1\times\sigma}=
V_{\partial_1\Delta^1\times\sigma}=\emptyset$ 
for all $\sigma\in\partial\Delta^{m+1}$.
\endroster
Hence $\{W_{\Delta^1\times\sigma}\to V_{\Delta^1\times\sigma}\}$ 
can
be considered as the
adic-surgery problem bounding the adic-surgery problem $\{\partial_0 W_{\Delta^1\times\sigma}\to \partial_0V_{\Delta^1\times\sigma}\}$.
The latter one represents $\partial_*[y]=[x]\in H_n(K,\LL)$.

Let $(f,b):W\to V$ be
the colimit understood as in $\S 3$. From Prop. 6.6 in \cite{LaMc}
we obtain a degree $1$
normal map $(f,b)$ between $(n+1)$-manifolds with boundaries
$$
\partial V=\partial_0 V\cup\partial' V \text{,  } \partial W=\partial_0 W\cup\partial' W
$$
where $\partial_0V$ is the colimit of $x$, $\partial' V$ is the colimit of the
restrictions to all $\partial_{m-|\sigma|+3}V_{\Delta^1\times\sigma}$, and similary for $\partial W$. 
Clearly, 
$f$ restricted to
$\partial' W$ is a homotopy equivalence. The colimit of the 
reference maps gives a reference map $(V,\partial_0V)\buildrel\SS\over\longrightarrow(B,K)$.

We now
describe
the map $p_*:H_n(K,\LL)\to H_n(B,\LL)$. Suppose that $[x]$ is 
represented by the normal degree $1$ map $(f,b):M\to X$ with the
reference map $\SS:X\to K$.
Then the family $M_{\sigma}\buildrel{f_\sigma,b_\sigma}\over\longrightarrow X_\sigma$, $\sigma\in K$, with $X_\sigma=\SS^{-1}(D(\sigma,K))$, $M_\sigma=f^{-1}(X_\sigma)$, $f_\sigma=f|_{M_\sigma}$ defines a cycle $x$. Note that
$p_*[x]$ is represented by the family $\widehat{M}_{\tau}\buildrel{\widehat{f}_\tau,\widehat{b}_\tau}\over\longrightarrow \widehat{X}_\tau$, $\tau\in B$, with $\widehat{X}_\tau=\widehat{\SS}^{-1}(D(\tau,B))$, $\widehat{M}_\tau=f^{-1}(\widehat{X}_\tau)$,  $\widehat{f}_\tau=f|_{\widehat{M}_\tau}$, where we have put $\widehat{\SS}=\SS\circ p$.

We have also to understand in this context the meaning of "homologies" between
cycles $x$, $x'$ on $K$, i.e. $[x]=[x']\in H_n(K,\LL)$.
There is a semi-simplicial map
$$
y:(\Delta^1\times\Sigma^{m},\Delta^1\times\overline{K})\to (\LL_{n-m},\emptyset)
$$
such that $y(\partial_0\Delta^1\times\sigma)=x(\sigma)$, $y(\partial_1\Delta^1\times\sigma)=x'(\sigma)$, i.e.
$y(\Delta^1\times\sigma)\in \LL_{n-m}(m-|\sigma|+1)$, defining surgery problems
$$
(F_{\Delta^1\times\sigma}, B_{\Delta^1\times\sigma}):W_{\Delta^1\times\sigma}\to V_{\Delta^1\times\sigma}
$$
of $(m-|\sigma|+5)$-ads as above, such that it restricts to
$$
(f_\sigma, b_\sigma):\partial_0W_{\Delta^1\times\sigma}=M_\sigma\to\partial_0V_{\Delta^1\times\sigma}=X_\sigma$$
$$
 \text{and } f'_\sigma, b'_\sigma:\partial_1W_{\Delta^1\times\sigma}=M'_\sigma\to\partial_1V_{\Delta^1\times\sigma}=X'_\sigma
$$
being $x(\sigma)$ and $x'(\sigma)$,
respectively. 
Moreover, $F_{\Delta^1\times\sigma}$ 
restricts to a homotopy equivalence
$$
\partial_{m-|\sigma|+3}W_{\Delta^1\times\sigma}\to \partial_{m-|\sigma|+3}V_{\Delta^1\times\sigma}.
$$
The colimit construction gives a normal cobordism $(F,B):W^{n+1}\to V^{n+1}$ 
such that $\partial W=\partial_0W\cup\partial_1W\cup\partial 'W \to\partial_0V\cup\partial_1V\cup\partial 'V$ with $\partial 'W\to\partial 'V$ is
a controlled homotopy equivalence extending  
the homotopy  equivalences $\partial M\to\partial X$ and $\partial M'\to\partial X'$.
Moreover,
$\partial_0W\to\partial_0V$
and
$\partial_1W\to\partial_1V$ represent $x$ (resp. $x'$) by normal degree $1$ maps.

We shall now consider the case of a closed topological $n$-manifold $K$. 
By the
Poincar\'{e} duality (cf. \cite{Ran}), we have
$$
H_n(K,\LL )=H^0(K,\LL )=[K,\Bbb Z\times G/TOP].
$$
The $\Bbb Z$-sector is related to Quinn's index invariant. 
The simply connected cover $\LL^+\to\LL$ is characterized 
by  $\LL^+_0=G/TOP$, and induces an injection $[K,G/TOP]\to [K,\Bbb Z\times G/TOP]$ 
onto the $1$-sector.
By duality it is the image of $H_n(K,\LL^+)\to H_n(K,\LL )$ 
(cf. \cite{Qu2}, \cite{Ran, $\S 25$}, and \cite{BFMW} for more details).

Given an element $[x]\in H_n(K,\LL^+)$ it defines a normal degree $1$ map 
between closed manifolds
$$
(g,c):Z^n\to K^n\text{ giving a cycle representation}$$
$$
\{z_\sigma\}=\{g_\sigma,c_\sigma\}:Z_\sigma\to K_\sigma=D(\sigma ,K)\subset D(\sigma,\partial\Delta^{m+1})$$

Obviously, $\partial_{m-|\sigma|+1}K_\sigma =\emptyset$, and the colimit of
the family $\{K_\sigma |\sigma\in K\}$ gives $K$. The cycle  $\{z_\sigma\}=z$ represents $[z]=[x]\in H_n(K,\LL^+)\subset H_n(K,\LL )$. 
If $(f,b):M^n\to X^n$ is the normal degree $1$ map obtained from the colimit of the cycle $x$,
there is a normal cobordism $C^{n+1}\to Y^{n+1}$ with $\partial_0C^{n+1}\to \partial_0Y^{n+1}$ 
equal to $(f,b)$ and $\partial_1C^{n+1}\to \partial_1Y^{n+1}$ equal to $(g,c)$. $\partial_0Y\cup \partial_1Y\subset\partial Y$, $\partial_0C\cup \partial_1C\subset\partial C$,
and $\partial Y\setminus(\partial_0Y\cup\partial_1Y)$ have
boundaries
$\partial 'X=\partial X$ 
defined above. Similarly, $\partial C\setminus(\partial_0C\cup\partial_1C)$ 
have
boundaries $\partial 'M=\partial M$.

We define $H_{n+1}^+(B,K,\LL )\subset H_{n+1}(B,K,\LL )$ as the
inverse image of $H_n(K,\LL^+)\subset H_n(K,\LL )$ 
under  $\partial_*:H_{n+1}(B,K,\LL )\to H_{n}(K,\LL )$. 
Obviously, $H_{n+1}(B,K,\LL^+)\subset H_{n+1}^+(B,K,\LL )$.

Finally, let  a class $[y]\in H_{n+1}^+(B,K,\LL )$ be given
with boundary cycle $x$, i.e. $[x]=\partial_*[y]$. Let
$$
W\to V\text{, } \partial_0W=M\buildrel{(f,b)}\over\longrightarrow\partial_0V=X $$
be the colimits of
$y$ and $x$. We glue $C\to Y$ to $W\to V$ along $M\buildrel{(f,b)}\over\longrightarrow X$ 
and obtain a normal degree $1$ map.
$$
(F^+,B^+):W^+\to V^+ \text{ with restriction} $$
$$
(g,c):\partial_0W^+=Z\to \partial_0V^+=K \text{ , } \partial 'W^+=\partial 'W \text{ , } \partial 'V^+=\partial 'V,$$
$$
(F^+,B^+)\text{ defines an element } [y^+]\in H_{n+1}^+(B,K,\LL ).$$

By the normal cobordism extension property (cf. \cite{Med, p. 45} or \cite{Wal, p. 93} 
applied in the proof of Theorem 9.6),
$W^+\to V^+$ and $W\to V$ are normally cobordant
by a cobordism being a product outside a small neighborhood of $C\to Y$.

In summary, we got for any closed oriented $n$-manifold $p:K\to B$ a
continuous 
map into a simplicial complex and any $[y]\in H_{n+1}^+(B,K,\LL )$ a cycle 
representation $y$ with colimit a degree $1$
normal map
$$
(f,b):W^{n+1}\to V^{n+1} $$
between $(n+1)$-dimensional manifolds such that $K=\partial_0V\subset\partial V$, $f$
restricted to $\partial_0W=M\to K$ is the colimit of a cycle of $\partial_*[y]\in H_n(K,\LL )$.
Moreover,
$f$ restricted to $\partial W\setminus\partial_0W\to \partial V\setminus\partial_0V$ 
is a homotopy equivalence. It
is a $\delta$-homotopy equivalence with respect to the reference map $\partial V\setminus\partial_0V\to B$. 

A special case is given when $[y]\in H_{n+1}^+(B,K,\LL )$ is in the image of $i_*:H_{n+1}(B,\LL )\to H_{n+1}^+(B,K,\LL )$, i.e. $i_*([z])=[y]$.
In this case one gets
a normal degree $1$ map $W\to V$ with $\partial_0W\to \partial_0V=\emptyset$ and $\partial W=\partial 'W\to \partial 'V=\partial V$ a controlled homotopy equivalence.

For our
later use we substitute $W\to V$ by adding $C\to Y$ representing the trivial 
element in $H_{n+1}^+(B,K,\LL )$. We can take $C=K\times I=Y$, with $\partial_0C=K\times 0=\partial_0Y$, $\partial_1C=\partial_1Y=\emptyset$ and $\partial 'C=K\times 1=\partial 'Y$.
We shall later refer to this case as the
"absolute case" and write it as
$$
W\dot{\cup}C\to V\dot{\cup}Y. $$
Moreover, we shall ignore the boundaries 
$\partial 'W$, $\partial 'V$, $\partial 'C$ and $\partial 'Y$.

\head {$\S$ 5. The map $a:H_{n+1}^+(B,K,\LL )\to \SS_{\varepsilon ,\delta}(K)$}\endhead
We assume $K$ to be a closed oriented $n$-manifold, $p:K\to B$ a $UV^1$-map, and 
$B$
a compact finite-dimensional metric ANR. To define the map $a$  it suffices to take
$B$ a finite ordered simplicial complex contained in $\partial\Delta^{m+1}$ as before.

We must first introduce the controlled structure set. Its definition
comes with the proof of the controlled surgery sequence of \cite{PQR}  (cf. also \cite{Fer}):
$$H_{n+1} (B,\LL )\to \SS_{\varepsilon ,\delta}(K \buildrel{p}\over\longrightarrow B) \to [K, G_{/TOP}] \to H_n (B,\LL ).$$
\indent The precise statement is: There exists $\varepsilon_0 >0$ such that for any $0<\varepsilon <\varepsilon_0$ there is  $\delta >0$ such that the sequence is exact provided $p$ is $a$ $UV^1$-map ($UV^1$($\delta$) would be sufficient). With these $\varepsilon$ and
$\delta$, the controlled structure set $\SS_{\varepsilon ,\delta}(K\to B)$ can be defined: Elements in $\SS_{\varepsilon ,\delta}(K\to B)$ are represented by pairs $(M, f)$, where $M$ is a closed $n$-manifold and $f:M\to K$ is a
$\delta$-homotopy equivalence with respect to the control map $p:K\to B$. 
Two pairs $(M_1, f_1)$ and  $(M_2, f_2)$ are "$\varepsilon$-related" if there is a homeomorphism $h:M_1\to M_2$ such that  $f_2\circ h$ and $f_1$ are $\varepsilon$-homotopic over $p:K\to B$ (i.e. the homotopy has "tracks" of size $<\varepsilon$ in $B$). Now, with this $\varepsilon$
and $f$, "$\varepsilon$-related" is an equivalence relation. We shall shortly write $\SS_{\varepsilon ,\delta}(K \buildrel{p}\over\longrightarrow B) = \SS_{\varepsilon ,\delta}(K)$.

\it{Further remarks:} \rm In \cite{PQR}, $\LL$ denotes the algebraic $\LL$-spectrum defined by adic quadratic Poincar\'{e} complexes (see \cite{Ran}). To a normal degree 1-map $(f,b):M\to K$ there is an obstruction in a controlled Wall
group $L_n (B,\varepsilon, \delta)$. It consists of quadratic $n$-dimensional Poincar\'{e} complexes of radius $\delta>0$ over $B$, modulo $\varepsilon$-cobordism. If this obstruction vanishes, controlled surgery can be completed in the middle dimension. The controlled Hurewicz and Whitehead theorems give an element in $\SS_{\varepsilon ,\delta}(K)$. All this involves $\varepsilon$-$\delta$-estimates. There is an assembly map $A:H_n (B, \LL)\to L_n(B,\varepsilon,\delta)$ defined by gluing the adic parts of an element in $H_n (B,\LL)$, similar as it was done in $\S3$. This map is proved to be an isomorphism for suitable $\varepsilon , \delta$. It is here, among other places, where the "stability threshold" $\varepsilon_0 >0$ drops in (see also \cite{RanYa}).

Let $[y]\in H_{n+1}^+(B,K,\LL )$ be given, and let $y$ be a cycle with $y(\sigma )$, given by 
normal degree $1$ $(m-|\sigma|+5)$-ads $(f_\sigma,b_\sigma):W_\sigma\to V_\sigma$.
Let $(f,b):W\to V$ be the colimit obtained in $\S 4$, i.e. $f$ restricts to a normal 
degree $1$ map between closed manifolds
$\partial_0W\to K$, and to a $\delta '$-homotopy equivalence $\partial 'W\to\partial 'V$. 
We claim that $\delta '$ can be choosen smaller than
$\delta$: First note that the above $\varepsilon_0>0$ depends on $B$ (and on $n$), 
and $\delta '$ becomes smaller and smaller if we subdivide
$B$ finer and finer. 

Let $\SS : V\to B$ be the reference map. 
Restricting $\SS$ to $\partial_0 V=K$ gives $p:K\to B$, and by assumption it is
$UV^1$.

We can assume that $S:V\to B$ is also $UV^1$:
First, change $W\to V$ as in the proof of Theorem 9.4
in \cite{Wal} to obtain an isomorphism 
$S_*:\pi_1(V)\to \pi_1(B)$.
Then since dim$V \ge 5$ one can approximate
$S$ (rel. $\partial V$)
to become a $UV^1$-map (cf. \cite{BFMW, Theorem 4.4}).

We can invoke the controlled simply-connected $\pi -\pi$-theorem (cf. e.g. \cite{Fer}):
We have a $\delta '$-equivalence $\partial_0W\to\partial_0V=K$, and we can 
assume that $\partial 'W\to\partial 'V$ is a $\delta '$-equivalence
too, and $\delta '$ can be as small as we need it. By the $\pi -\pi$-theorem
there are $\delta_0 >0$, $k>0$, such that for $\delta '<\delta_0$
there is a normal cobordism between $(W,\partial_0W)\to (V,\partial_0V)$ 
relative $\partial 'V$ to a $k\delta '$-homotopy equivalence
$$
f':(W',\partial_0W')\to (V,\partial_0V=K).$$
Choosing further $k\delta '<\delta$, we get a class in $\SS_{\varepsilon , \delta}(K)$ 
represented by
$$
f_0':\partial_0W'\to K\text{ , }f_0'=f'|_{\partial_0W'}.$$
We shall define $a:H_{n+1}^+(B,K,\LL )\to \SS_{\varepsilon , \delta}(K)$ 
by $a([y])=[f_0',\partial_0W']$.

It remains to show that $a$ is well-defined. Suppose $y$, $y_1$ are cycles with $[y_1]=[y]$, let
$$
(W,\partial_0W)\buildrel{(f,b)}\over\longrightarrow (V,\partial_0V=K)$$
$$
(W_1,\partial_0W_1)\buildrel{(f_1,b_1)}\over\longrightarrow (V_1,\partial_0V_1=K)$$
be colimits 
and let 
$$
(W',\partial_0W')\buildrel{f'}\over\longrightarrow (V,\partial_0V=K)$$
$$
(W_1',\partial_0W_1')\buildrel{f'_1}\over\longrightarrow (V_1,\partial_0V_1=K)$$
be $k\delta '$-, resp. $k\delta_1 '$-homotopy equivalences obtained from the
controlled $\pi -\pi$-theorem. We may assume 
$k\delta '<\delta$, $k\delta_1 '<\delta$.

Let $Q\to \Omega$ be a normal cobordism between $(f,b)$ and $(f_1,b_1)$. It
restricts to a normal cobordism $P\to K\times I$
with $\partial P=M\cup M_1\buildrel{f\cup f_1}\over\longrightarrow K\times\{0,1\}$,
and with boundary 
$\partial Q=W\cup P\cup W_1\to \partial\Omega =V\cup K\times I\cup V_1$.

To proceed we briefly have to recall the proof of the
$\pi -\pi$-theorem: 
The first step is to make $(W,M=\partial_0W)\to (V,\partial_0V=K)$ highly connected,
which we can assume.

The second step consists of handle
subtraction in the pairs $(W,M)$ and $(W_1,M_1)$,
i.e. $(W',\partial_0W')$ is obtained from $W$ by
subtracting handles. We add these handles to $P$ along $M$. Similarly,
we add the 
handles subtracted from $(W_1,M_1)$ along $M_1$ to $P_1$.
One obtains $P'$ with $\partial P'=M'\cup M_1'$. In other words,
we obtain a new
decomposition of $\partial Q$ and $\partial\Omega$ and
$\partial Q\to\partial\Omega$ as
$$
\partial Q=W'\cup P'\cup W_1'\to \partial\Omega =V\cup K\times I\cup V_1.$$

Hence $Q\to\Omega$ displays a normal bordism between 
$$
(P',\partial P'=M'\cup M_1')\to (K\times I, K\times\{0\}\cup K\times\{1\})$$
and the controlled homotopy equivalence
$$
(W',M')\cup (W_1',M_1')\to (V,K)\cup (V_1,K).$$

We denote by
$\phi ':P'\to K\times I$ the first map.
Then $\phi '$
restricted to $M'$, $M_1'$ is a $k\delta '$ (resp. $k\delta_1'$) equivalence,
with $k\delta '$, $k\delta_1'<\delta$.

Let us assume $\delta_1'=\delta '$.
Controlled surgery theory then implies (cf. \cite{PQR}, \cite{Qu1}, \cite{Fer})
that $\phi '$ can be surgered to a simple $\overline{k}k_1\delta '$-homotopy equivalence
$$
H:(P'',M',M'_1)\to (K\times I, K\times\{0\},  K\times\{1\})$$
(for some $\overline{k}>0$) with $H|_{M'}=f'$, $H|_{M_1'}=f_1'$.
In particular $P''$ is an $s$-cobordism.

Let $G:P''\to M'\times [0,1]$ be a homeomorphism with
$$
G|_{M'}=Id \text{ , } g=G|_{M_1'} :M_1'\to M. $$

The composition
$$
M_1'\times I \buildrel{g\times Id}\over\longrightarrow M'\times I\buildrel{G^{-1}}\over\longrightarrow P''\buildrel{H}\over\longrightarrow K\times I
$$
restricts to $f'\circ g:M_1'\times 0\to K\times 0$ and $f_1':M_1'\times \{1\}\to K\times \{1\}$. 
Therefore $f'\circ g$ and $f_1'$ are
$\overline{k}k\delta '$-homotopic over $B$.

However,
we can choose $\delta '$ such that $\overline{k}k\delta '<\varepsilon$, 
showing $f':M'\to K$, $f_1':M_1'\to K$ define the same element in 
$\SS_{\varepsilon , \delta }(K)$, hence the map $a$ is well-defined.

This finishes the generic case, and we shall now come to the absolute case.
Now let $[y]\in H_{n+1}^+(B,K,\LL )$ be the image of $[z]\in H_{n+1}^+(B,\LL )$, i.e. $i_*[z]=[y]$. 
Then we get a normal degree $1$ map
$W\dot{\cup}C\to V\dot{\cup}Y$ 
such that $\partial_0V=\emptyset$, $\partial_0W=\emptyset$, $C\to Y$ 
is a controlled homotopy equivalence with $\partial C=K\to \partial Y=K$ a homeomorphism.
The $\pi -\pi$-theorem cannot be applied to $(W,\partial_0W)$ since
$\partial_0W=\emptyset$. 
We define $a([y])$ as the
result of the controlled Wall action of $[z]$ on $Id:K\to K$, 
i.e. $a[y]$ is the image of $[z]$ under
$$
H_{n+1}(B,\LL )\to \SS_{\varepsilon, \delta}(K)$$
(cf. \cite{PQR}).
We have to show that the map $a$ is well-defined.
It suffices
to consider the following
two types 
of cycle representations of $i_*([z])=[y]$:
\roster
\item "(I)"
$ W\dot{\cup}C\to V\dot{\cup}Y$  is an absolute representation;
\item "(II)"
$ W_1\to V_1$ is a generic representation;
\endroster
with $\partial_0W_1\to\partial_0V_1\not=\emptyset$.

Let $Q\to \Omega$ be a normal cobordism between them.
Then $\partial Q=W\dot{\cup}C\mathop\cup\limits_{K}P\cup W_1\to 
\partial\Omega = V\dot{\cup}Y\mathop\cup\limits_{K} K\times I\cup V_1,$ where $(P,\partial_0P=K,\partial_1P=\partial_0W_1)\to (K\times I, K\times 0, K\times 1)$.

The restriction to $W_2=C\cup P\cup W_1\to V_2=Y\cup K\times I\cup V_1$ is a normal 
degree $1$ map with controlled homotopy equivalence on the boundary. It defines the
same element $i_*[z]=[y]\in H_{n+1}^+(B,K,\LL )$. We may also cancel $C\to Y$ thus
obtaining 
$$
W_3=P\cup W_1\to V_3\cong V_1 \text{, and identity on}$$
$$
\partial_0W_3=K\to \partial_0V_3=K.$$

It still defines the same element $[y]$. We can apply the $\pi -\pi$-theorem to $(W_3,\partial_0W_3=K)\to (V_3,\partial_0V_3=K)$ and obtain a controlled homotopy 
equivalence
$$
(W_3',\partial_0W_3')\to (V_3,\partial_0V_3=K).$$

We have to show that the $\delta$-homotopy equivalence $\partial_0W_3'\to K$ defines
the same element as the image of
$$
[z]\in H_{n+1}(B,\LL )\to \SS_{\varepsilon ,\delta}(K).$$

First one makes $W_3\to V_3$ highly connected by surgeries in the interior of $W_3$. We 
shall assume this for $W_3\to V_3$. Assume $n+1=2k$.
The middle-dimensional classes defining the surgery kernel can be represented by 
framed immersed
$(k+1)$-handles
$$
(H', \partial H')\to (W_3, K).$$

Let $U=K\times [0,\rho]\cup H' \subset W_3$, where $K\times [0,\rho]$ denotes a small 
collar of $K$. Note
$U$ can be identified with the controlled Wall action of $[z]$ on $K$. On the other 
hand one can remove the intersection of the immersed handles across $K$ (here $\pi_1(K)\mathop\to\limits_{\cong} \pi_1W_3$ is used) to get embedded handles $(H,\partial H)\subset (W_3, K)$. 

We may assume that $U=K\times [o,\rho]\cup H$.
Then $(W_3', \partial_0W_3')=(W_3\setminus \buildrel{\circ}\over{U}, \partial_1U)\to (V_3,K)$ 
is the result from the $\pi -\pi$-theorem. This shows that both definitions
of the map $a$ on $i_*[z]=[y]$ 
give the same result.

The case when $n+1$ is odd goes as follows:
One has to consider neighbourhoods $U$ of $K$ 
in $W_3$ given by attached embedded handles of the
type $S^k\times D^{k+1}$, $2k+1=n+1$. 
They have to be subtracted. On the other hand,
in order
to realize an element of $H_{n+1}(B,\LL )$
one attaches trivial handles and then
applies the controlled obstruction to the 
transverse spheres $1\times S^k$. On this result one then performs the surgeries.
The resulting homotopy equivalence is the same.

\head {$\S$ 6. Proof of Theorem 3}\endhead

The map $b$ was
defined in \cite{DrFe} as the
composite map
$$
H_{n+1}(B, K,\LL )\cong S_{n+1}(B,K,\LL )\buildrel{\partial_*}\over\longrightarrow S_n(K,\LL )\to\SS (K).
$$
Here, $S_{n+1}(B,K,\LL )$, $S_n(K,\LL )$ are the structure sets defined in \cite{Ran} (however,
denoted differently as $S_{n+2}(B,K)$, $S_{n+1}(K)$).

An element in $S_n(K,\LL)$ is represented by a homotopy equivalence of $n$-manifolds $P^n\to Q^n$
with the
reference map $Q^n\to K$. There is a natural inclusion $\SS (K)\subset S_n(K,\LL )$.
The first isomorphism is due to the $\pi -\pi$-theorem. The map $S_n(K,\LL )\to \SS ( K)$ 
is such that the composite $\SS (K)\subset S_n(K,\LL )\to \SS (K)$ is the identity.

Let $[y]\in H_{n+1}^+(B,K,\LL )$. Our construction of $a[y]$ produces a homotopy equivalence $(W',\partial_0 W')\to (V,\partial_0V=K)$ with 
the reference map $(V,\partial_0V=K)\to (B,K)$. 

This is the image of $[y]$ in $S_{n+1}(B,K,\LL )$, hence $\partial_0W'\to \partial_0V$ 
is the image in $S_n(K,\LL )$.
However,
it belongs already to $\SS(K)\subset S_n(K,\LL )$, and $\SS(K)\subset S_n(K,\LL )\to \SS(K)$
is the identity. This shows that the following diagram
$$
\xymatrix @C=.2cm @R=.1cm{
H_{n+1}(B, K,\LL ) &\cong & S_{n+1}(B,K,\LL ) & \to  & S_n(K,\LL )  & \to & \SS(K)  \\
\cup \\
H_{n+1}^+(B, K,\LL ) \ar[rrrrrr]^{a}  & & & & &  & \SS_{\varepsilon ,\delta}(K\to B) \ar[uu]
}
$$
commutes.

To complete the proof  of Theorem 3 we have to verify the
commutativity of the diagrams
$$
\xymatrix @R=.2cm{
& H_{n+1}(B,\LL) \ar@{=}[dd] \ar[r] & H_{n+1}^+(B, K,\LL ) \ar[dd]^a \\
(I) & & \\
& H_{n+1}(B,\LL) \ar[r] & \SS_{\varepsilon,\delta}(K\to B) 
}$$
and
$$
\xymatrix @R=.2cm{
& H_{n+1}^+(B,K, \LL) \ar[dd]^{a} \ar[r] & H_{n}(K,\LL^+) \ar[dd]^{\cong} \\
(II) & & \\
& \SS_{\varepsilon,\delta}(K\to B)  \ar[r] & [K,G/TOP] .
}$$

Now, (I) follows from the definition. To prove commutativity of (II), we consider the following big diagram

$$
\xymatrix @C=.2cm @R=.1cm{
S_{n+1}(B,K,\LL ) \ar[ddd]^{\cong} &  \ar[rrrrrrrrrr]^{\partial_*} & & & & & & & & & & & S_{n}(K,\LL ) \ar[ddd]\\
& & & & & & & & & & & &\\
& & & & & & & & & & & &\\
H_{n+1}(B,K,\LL ) \ar@{=}[ddd] &  \ar[rrrrrrrrrr]^{\partial_*} & & & & & & & & & & & H_{n}(K,\LL ) \ar@{=}[ddd]\\
& & & & & & & & & & & &\\
& & & & & & & & & & & &\\
H_{n+1}(B,K,\LL )  \ar[ddd]^{b} & & \supset & & H_{n+1}^+ (B,K,\LL ) \ar[ddd]^{a} &  \ar[rr]^{\partial_*} & & & H_n (K, \LL^+) \ar[ddd]^{\cong} & \ar[rr] & & & H_n (K, \LL) \ar[ddd]^{\cong} \\
& & & & & & & & & & & &\\
& & & & & & & & & & & &\\
S(K)  \ar@{=}[ddd] & & & \ar[ll] & S_{\varepsilon, \delta}(K) & \ar[rr] & & & [K, G_{/TOP}]  \ar@{=}[ddd]  & & &  \ar[ll]_{pr} & [K, Z\times G_{/TOP}]\\
& & & & & & & & & & & &\\
& & & & & & & & & & & &\\
S(K) &\ar[rrrrrr]^{\eta} & & & & & & & [K, G_{/TOP}] & & & &
}
$$

The "inner" diagram is (II). The upper 3 rows are just definitions of various $\partial_*$ (see also \cite{Ran}, p. 207). Commutativity of the middle
left
diagram is explained in the above diagram. The middle
right
diagram commutes because of compatibility of $\LL$ and $\LL^+$-Poincar\'{e} duality. The lower diagram commutes, here $\eta$ associates to a
homotopy equivalence $f:M\to K$ its normal invariant, and $pr$ denotes the projection map.
The composed maps of the "outer" diagram give the diagram

$$
\xymatrix @C=.2cm @R=.1cm{
S_{n+1}(B,K,\LL ) \ar[ddd] &  \ar[rrrr]^{\partial_*}  & & & & & S_{n}(K,\LL ) \ar[ddd]\\
& & & & &  &\\
& & & & & & \\
S(K) &\ar[rrrr]^{\eta}  & & & & & [K, G_{/TOP}].
}
$$

This commutes because of commutativity of

$$
\xymatrix @C=.2cm @R=.1cm{
H_{n}(K,\LL )  &  \ar[rrrr]^{\cong}  & & & & & [K, G_{/TOP}] \ar[ddd]^{pr}\\
& & & & &  &\\
& & & & & & \\
H_n(K, \LL^+) \ar[uuu] &\ar[rrrr]^{\cong}  & & & & & [K, G_{/TOP}].
}
$$

Commutativity of (II) follows from commutativity of all these diagrams.

\head{\smc Acknowledgements} \endhead
   This research was supported by the Slovenian Research Agency grants P1-0292-0101 and J1-4144-0101. We thank the referee for comments.

\enddocument